\documentclass{amsart}
\usepackage{graphicx}

\newtheorem{theorem}{Theorem}[section]
\newtheorem{lemma}[theorem]{Lemma}
\newtheorem{corollary}[theorem]{Corollary}

\theoremstyle{definition}

\newtheorem{example}[theorem]{Example}

\newtheorem{claim}[theorem]{Claim}

\theoremstyle{remark}

\numberwithin{equation}{section}

\begin{document}

\title[Rational structure on algebraic tangles]{Rational structure on algebraic tangles and closed incompressible surfaces in the complements of algebraically alternating knots and links}

\author{Makoto Ozawa}
\address{Department of Natural Sciences, Faculty of Arts and Sciences, Komazawa University, 1-23-1 Komazawa, Setagaya-ku, Tokyo, 154-8525, Japan}
\email{w3c@komazawa-u.ac.jp}





\begin{abstract}
Let $F$ be an incompressible, meridionally incompressible and not boundary-parallel surface with boundary in the complement of an algebraic tangle $(B,T)$. Then $F$ separates the strings of $T$ in $B$ and the boundary slope of $F$ is uniquely determined by $(B,T)$ and hence we can define the slope of the algebraic tangle. In addition to the Conway's tangle sum, we define a natural product of two tangles. The slopes and binary operation on algebraic tangles lead an algebraic structure which is isomorphic to the rational numbers.

We introduce a new knot and link class, algebraically alternating knots and links, roughly speaking which are constructed from alternating knots and links by replacing some crossings with algebraic tangles. We give a necessary and sufficient condition for a closed surface to be incompressible and meridionally incompressible in the complement of an algebraically alternating knot or link $K$, in particular we show that if $K$ is a knot, then the complement of $K$ does not contain such a surface.
\end{abstract}

\maketitle

\section{Introduction}

Conway introduced rational tangles and algebraic tangles for enumerating knots and links by using his ``Conway notation'' (\cite{C}). He noted that the rational tangles correspond to the rational numbers in one-to-one.
Kauffman and Lambropoulou gave also a proof of this theorem (\cite{KL}).
Gabai gave an another definition for algebraic links by a plumbing construction from a weighted tree, so these links are called an ``arborescent'' links (\cite{G}).

Hatcher and Thurston classified incompressible surfaces in the complements of 2-bridge knots (\cite{HT}).
Oertel classified closed incompressible surfaces in the complements of Montesinos links (\cite{O}).
Bonahon and Siebenmann characterized non-hyperbolic algebraic links (\cite{BS}).
Wu classified non-simple algebraic tangles (\cite{W}) and Reif determined the hyperbolicity of algebraic tangles and links (\cite{R}).
Futer and Gu\'{e}ritaud gave another proof of the theorem by Bonahon and Siebenmann (\cite{FG}).
In another direction, Menasco showed that closed incompressible surfaces in alternating link complements are meridionally compressible, and that an alternating link is split if and only if the alternating diagram is split (\cite{M}).

Krebes discoverd that the greatest common divisor of the determinant of a numerator and denominator of a tangle embedded in a link divides the determinant of the link (\cite{K}).
Moreover he constructed a map from tangles to formal fractions (not necessarily reduced), and the map on the algebraic tangles is surjective.

In this paper, we extend the slope of rational tangles to one of algebraic tangles by means of the boundary slope of essential surfaces in algebraic tangles, and show that a map from algebraic tangles to the boundary slopes coincides with ``reduced'' Krebes's invariant.
This map induces a homomorphism from the algebraic tangles to the rational number.
We also introduce algebraically alternating knots and links which are defined by means of diagrams with Conway notation of ``alternating slope'', and show that any essential closed surface in the complement of algebraically alternating knots is meridionally compressible.

\subsection{Rational structure on algebraic tangles}

Let $M$ be a 3-manifold and $T$ a 1-manifold properly embedded in $M$.
We say that a surface $F$ properly embedded in $M-T$ is {\em m-essential} (meridionally essential) if it is incompressible, meridionally incompressible and not boundary-parallel in $M-T$.
For the definition of incompressible and meridionally incompressible, we refer to \cite{MO}.
(In general, the condition that it is boundary-incompressible is required, but under the situation in this paper, it follows from the above three conditions.
Also surfaces are automatically orientable under the situation of this paper.)





\medskip
\noindent{\bf Theorem \ref{separating}}
{\em Let $(B,T)$ be an algebraic tangle.
Then, any m-essential surface with boundary in $B-T$ separates the components of $T$, and all boundary slopes of m-essential surfaces with boundary in $B-T$ are unique.}
\medskip

In the next section, we will see that there exists at least one m-essential surface with boundary in any algebraic tangle (Lemma \ref{sum}).
Hence we can define the {\em slope} of an algebraic tangle $(B,T)$ as the boundary slope of an m-essential surface $F$ in $B-T$.
This definition succeeds the slope of rational tangles.
We will state the definitions of the summation and multiplication for tangles in the next section.

\medskip
\noindent{\bf Theorem \ref{homomorphism}}
{\em Let $\phi$ be a map from the set of algebraic tangles to the set of rational numbers which maps an algebraic tangle $(B,T)$ to the slope of $(B,T)$.
Then $\phi$ is a surjective homomorphism under the summation and multiplication for tangles.}
\medskip

The homomorphism $\phi$ is consistent with the Krebes invariant $f$ (\cite{K}), to be precise, $\phi(T)$ is equal to an irreducible fraction of $f(T)$.

\subsection{Algebraically alternating knots and links}

\begin{figure}[htbp]
	\begin{center}
	\begin{tabular}{cc}
	\includegraphics[trim=0mm 0mm 0mm 0mm, width=.4\linewidth]{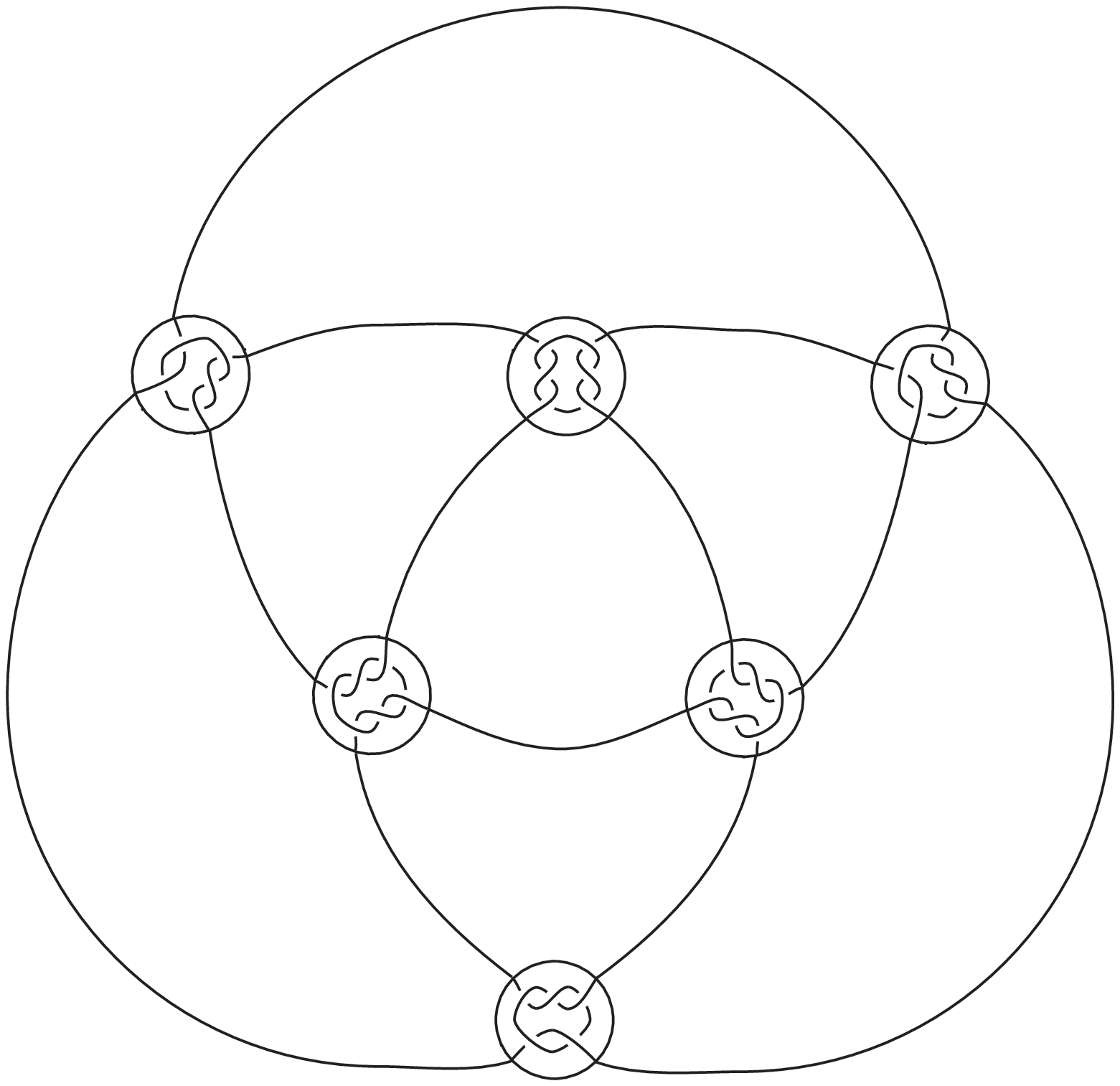}&
	\includegraphics[trim=0mm 0mm 0mm 0mm, width=.4\linewidth]{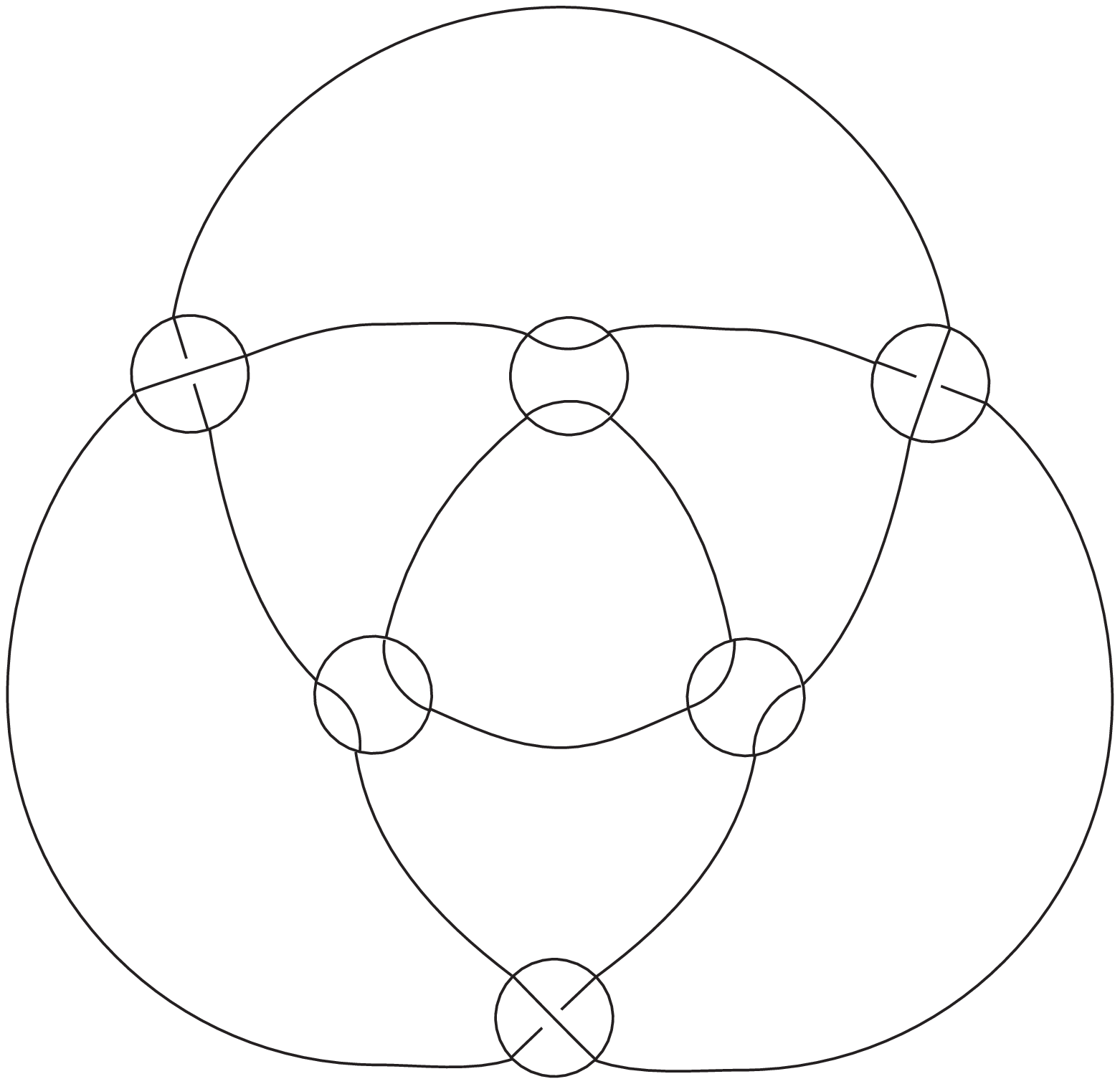}\\
	$\tilde{K}$ & $\tilde{K_0}$
	\end{tabular}
	\end{center}
	\caption{an algebraically alternating link diagram $\tilde{K}$ and the basic diagram $\tilde{K_0}$}
	\label{algebraically}
\end{figure}

Let $S^2\subset S^3$ be the standard 2-sphere and $G$ a connected quadrivalent graph in $S^2$ which has no bigon.
We obtain a knot or link diagram $\tilde{K}$ by substituting algebraic tangles for vertices of $G$.
After such an operation, we substitute each algebraic tangle $(B,T)$ for a rational tangle of slope $1$, $-1$, $0$ or $\infty$ if the slope of $(B,T)$ is positive, negative, $0$ or $\infty$ respectively (fixing four points of $\partial T$).
The resultant knot or link diagram is said to be {\em basic} and denoted by $\tilde{K_0}$.
Then we say that $\tilde{K}$ is {\em algebraically alternating} if $\tilde{K_0}$ is alternating, and $K$ is {\em algebraically alternating} if $K$ has an algebraically alternating diagram.
See Figure \ref{algebraically}.
We remark that the class of algebraically alternating knots and links includes both of alternating knots and links and algebraic knots and links.

Besides the slope, we can define the {\em genus} of an algebraic tangle $(B,T)$ as the minimal genus of m-essential surfaces with boundary in $B-T$.

\begin{figure}[htbp]
	\begin{center}
	\includegraphics[trim=0mm 0mm 0mm 0mm, width=.6\linewidth]{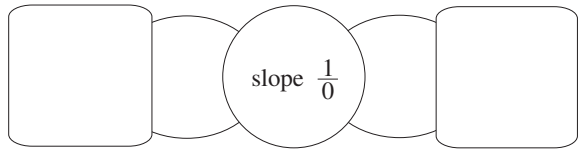}
	\end{center}
	\caption{Cut tangle}
	\label{cut}
\end{figure}

In a diagram of Conway notation, an algebraic tangle of slope $\infty$ (resp. $0$) is called a {\em cut tangle} if the diagram becomes split when we replace the algebraic tangle with a rational tangle of slope $1/0$ (resp. $0/1$).
See Figure \ref{cut}.

\begin{figure}[htbp]
	\begin{center}
	\includegraphics[trim=0mm 0mm 0mm 0mm, width=.25\linewidth]{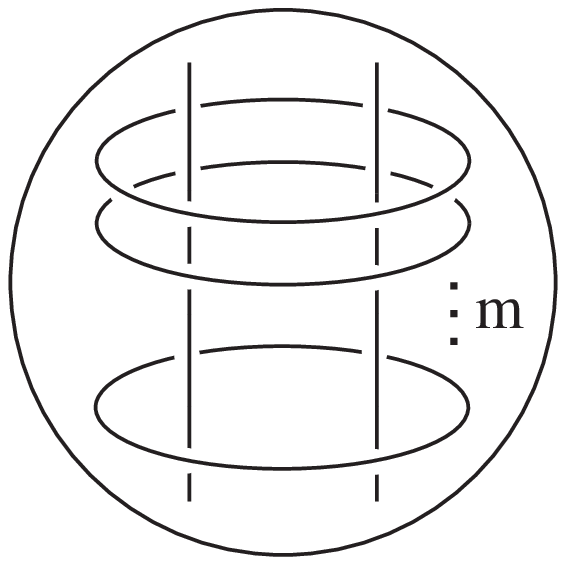}
	\end{center}
	\caption{$Q_m$}
	\label{Q_2}
\end{figure}

An algebraic tangle with two unknotted parallel strings and $m$ unknotted parallel loops is denoted by $Q_m$ (Figure \ref{Q_2}).

Menasco (\cite{M}) showed that there exists no m-essential closed surface in the complement of any alternating knot or link.
Oertel (\cite{O}) showed that the complement of a Montesinos knot or link $K$ containes an m-essential torus if and only if $K$ is a pretzel link $P(p,q,r,-1)$, where $\frac{1}{p}+\frac{1}{q}+\frac{1}{r}=1$.
Bonahon and Seibenmann (\cite{BS}) showed that the complement of a large algebraic link contains an m-essential torus if and only if $(S^3,K)$ contains $Q_2$.
Here, an algebraic link is {\em large} if it does not have a form of a Montesinos link with length 3.
See also \cite{FG} and \cite{R}.

Next theorem extends the above results.

\begin{theorem}\label{closed}
Let $K$ be an algebraically alternating knot or link in $S^3$, $\tilde{K}$ an algebraically alternating diagram of $K$ and $\tilde{K_0}$ the basic diagram of $\tilde{K}$.
Then,
\begin{enumerate}
\item There exists an m-essential closed surface in $S^3-K$ if and only if $\tilde{K_0}$ is split or there exists an algebraic tangle in $\tilde{K}$ which contains a closed algebraic sub-tangle.
\item There exists an m-essential 2-sphere in $S^3-K$ if and only if there exists a genus $0$ cut tangle in $\tilde{K}$.
\item Suppose that there exists no genus $0$ cut tangle in $\tilde{K}$.
Then there exists an m-essential torus in $S^3-K$ if and only if there exists a genus $1$ cut tangle in $\tilde{K}$ or $(S^3,K)$ contains $Q_2$.
\end{enumerate}
\end{theorem}

Next corollary extends the Menasco's meridian lemma (\cite{M}).

\begin{corollary}
Let $K$ be an algebraically alternating knot, and $F$ a closed incompressible surface in the complement of $K$.
Then, $F$ is meridionally compressible.
\end{corollary}

\begin{proof}
Suppose that $\tilde{K_0}$ is split.
Since $K$ consists of one component, any algebraic tangle of $\tilde{K}$ contains no loop component.
Therefore by Lemma \ref{connection}, $\tilde{K}$ has more than one component.
Next suppose that there exists an algebraic tangle in $\tilde{K}$ which contains a closed algebraic sub-tangle.
Then by Lemma \ref{connection}, some algebraic tangle of $\tilde{K}$ has a loop component.
In either cases, we have a contradiction.
\end{proof}

By \cite[Lemma 2.1]{IO}, we have the following corollary.

\begin{corollary}
Hyperbolic algebraically alternating knots satisfy the Menasco-Reid conjecture $($\cite{MR}$)$, i.e. there is no totally geodesic closed surface embedded in the knot complement.
\end{corollary}

\begin{figure}[htbp]
	\begin{center}
	\begin{tabular}{cc}
	\includegraphics[trim=0mm 0mm 0mm 0mm, width=.3\linewidth]{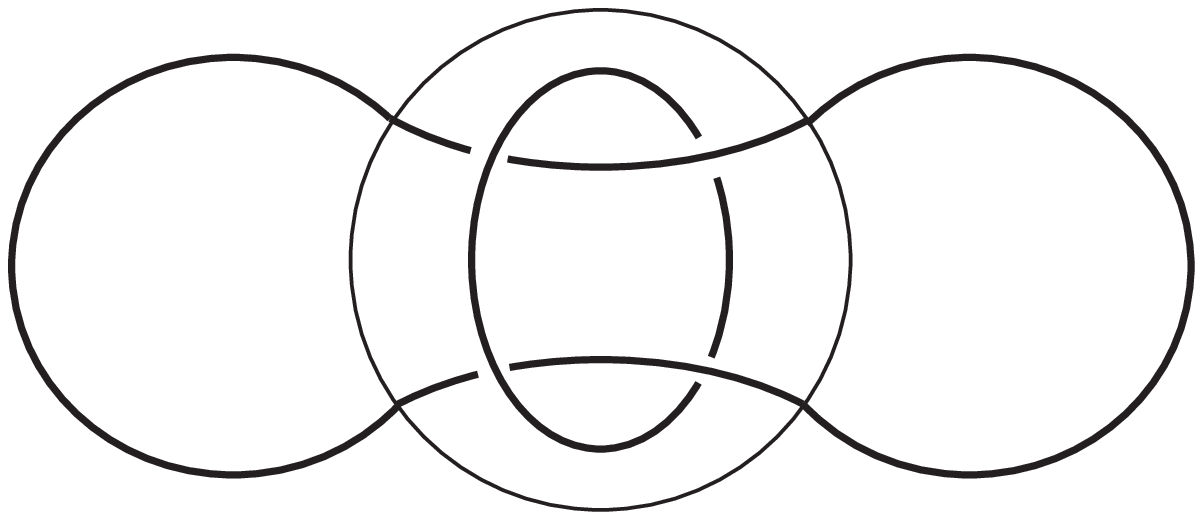}&
	\includegraphics[trim=0mm 0mm 0mm 0mm, width=.3\linewidth]{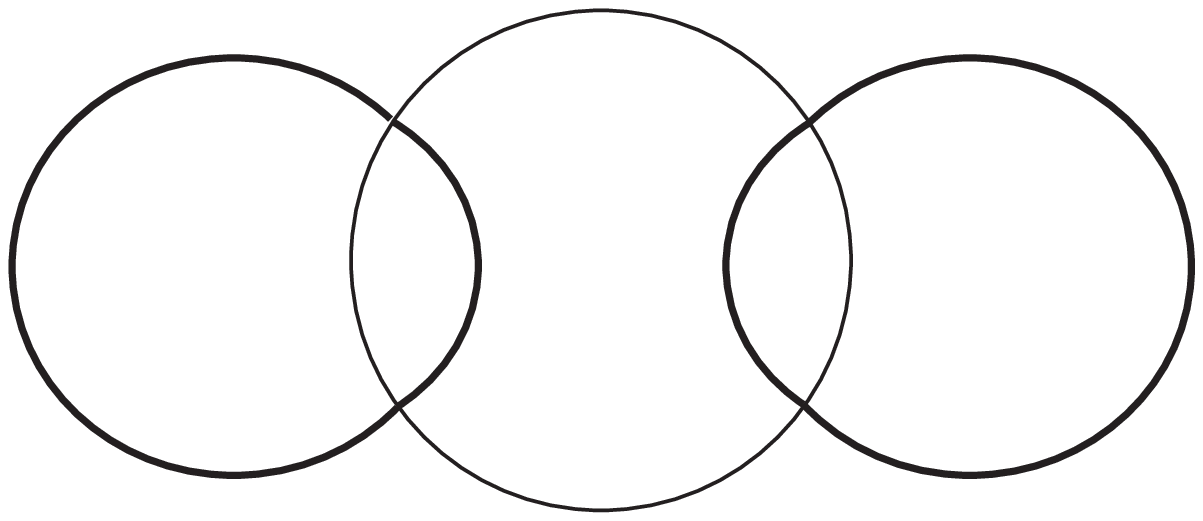}\\
	$\tilde{K}$ & $\tilde{K_0}$
	\end{tabular}
	\end{center}
	\caption{an algebraically alternating link diagram $\tilde{K}$ and the basic diagram $\tilde{K_0}$}
	\label{split}
\end{figure}

\begin{example}
Let $\tilde{K}$ be an algebraically alternating diagram of $K$ as Figure \ref{split}.
Then $\tilde{K}$ contains a genus 0 cut tangle since $S^3-K$ contains an m-essential 2-sphere.
\end{example}

\begin{figure}[htbp]
	\begin{center}
	\begin{tabular}{cc}
	\includegraphics[trim=0mm 0mm 0mm 0mm, width=.3\linewidth]{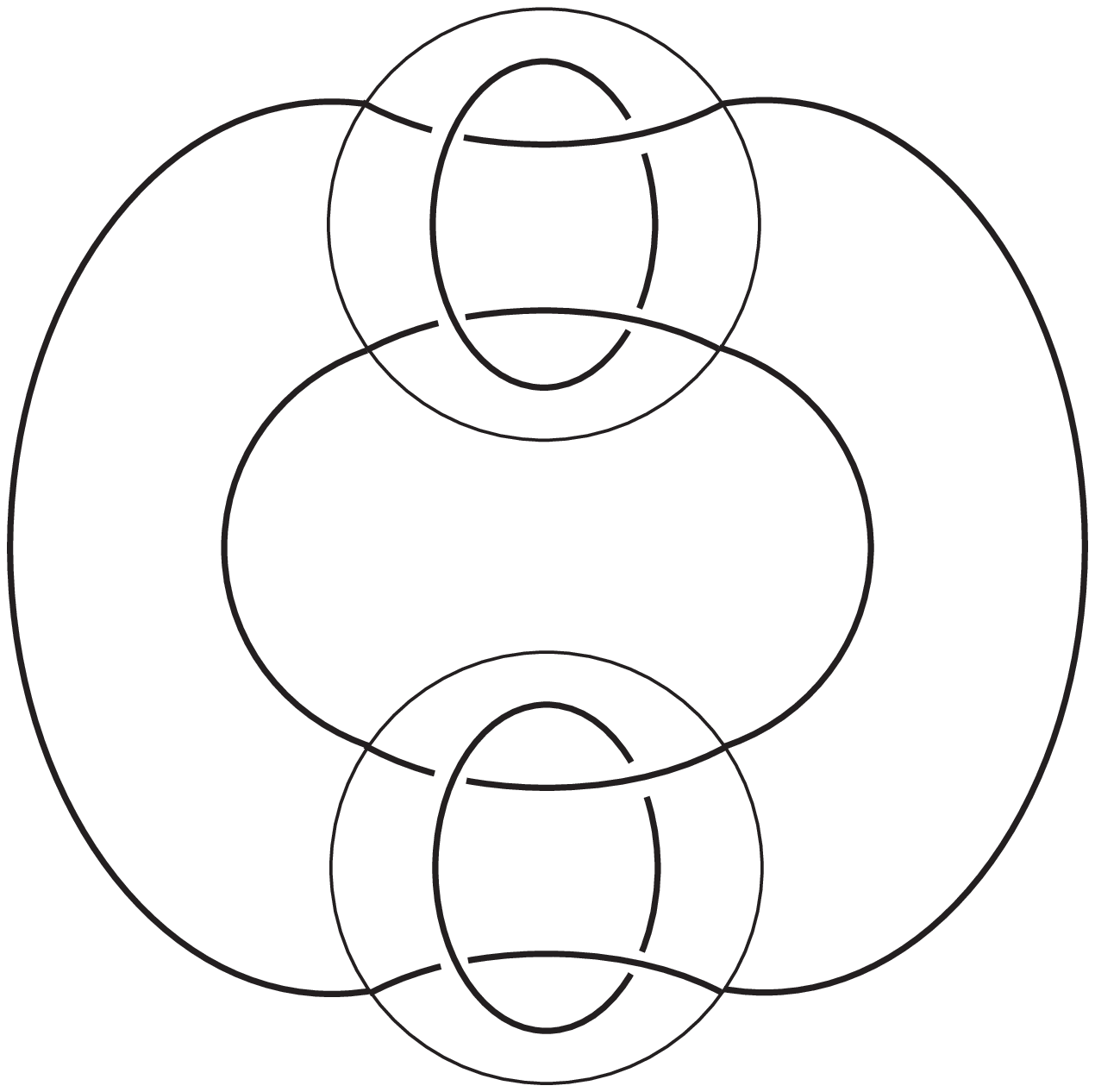}&
	\includegraphics[trim=0mm 0mm 0mm 0mm, width=.3\linewidth]{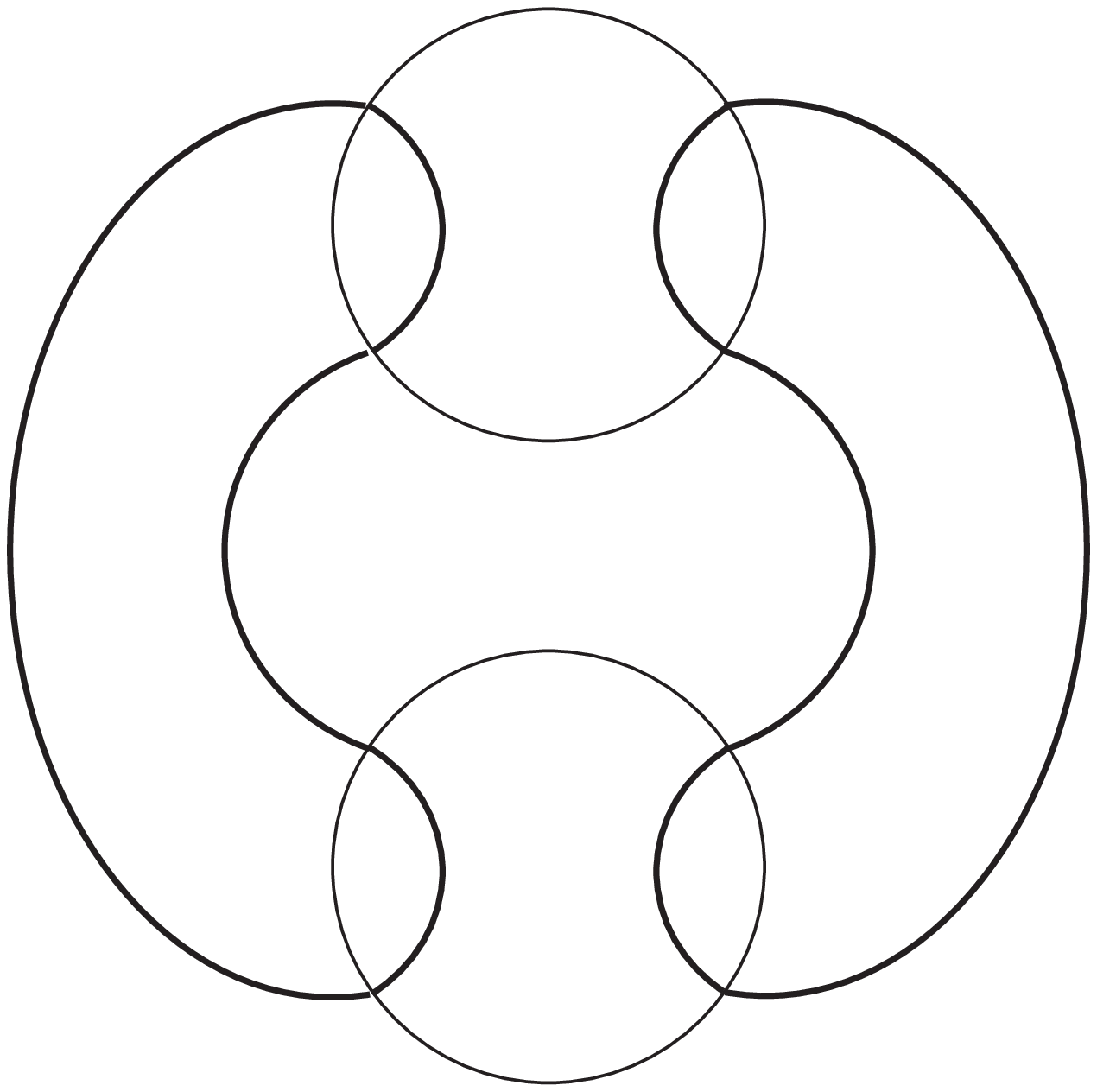}\\
	$\tilde{K}$ & $\tilde{K_0}$
	\end{tabular}
	\end{center}
	\caption{an algebraically alternating link diagram $\tilde{K}$ and the basic diagram $\tilde{K_0}$}
	\label{torus}
\end{figure}

\begin{example}
Let $\tilde{K}$ be an algebraically alternating diagram of $K$ as Figure \ref{torus}.
Then $(S^3,K)$ contains $Q_2$ since $S^3-K$ contains an m-essential torus and $\tilde{K}$ contains no genus 0 nor 1 cut tangle.
\end{example}

\begin{figure}[htbp]
	\begin{center}
		\includegraphics[trim=0mm 0mm 0mm 0mm, width=.7\linewidth]{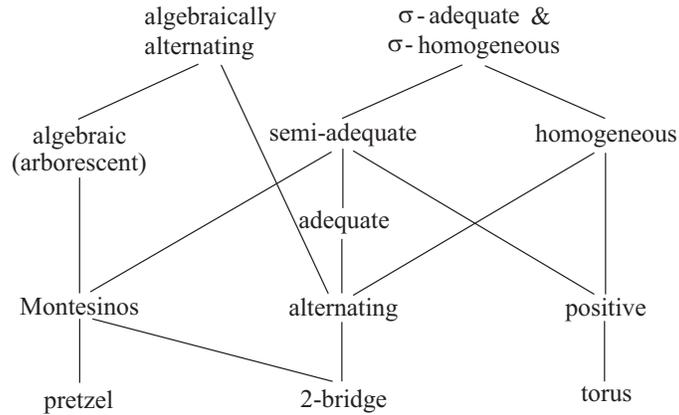}
	\end{center}
	\caption{The Hasse diagram of various classes of links}
	\label{benz}
\end{figure}

Here we represent the Hasse diagram of various classes of knots and links which are defined by means of diagrams.
See Figure \ref{benz} and refer Montesinos links (\cite{O}), alternating links (\cite{M}), positive links (\cite{MO1}), homogeneous links (\cite{C1}), adequate links (\cite{T}), semi-adequate links (\cite{LT}), algebraic links (\cite{C}), arborescent links (\cite{G}), $\sigma$-adequate and $\sigma$-homogeneous links (\cite{MO2}) for the definition of each diagrams.

\section{Rational structure on algebraic tangle}

\subsection{The slope of rational tangles}

Let $(B,T)$ be a 2-string tangle fixing four points of $\partial T$.
There are a double covering map $d:S^1\times S^1\to \partial B$ branched over $\partial T$ and a universal covering map $u:\Bbb{R}^2\to S^1\times S^1$.
Then we have a map $p:\Bbb{R}^2\to \partial B$ as a composition of these two covering maps such that $p^{-1}(\partial T)$ is the set of half integral points.
We say that an essential loop $C$ in $\partial B-\partial T$ which separates $\partial T$ into two pairs of two points has a {\em slope $p/q$} if a component of $p^{-1}(C)$ is a line with slope $p/q$.
In particular, $C$ is a {\em meridian} (resp. a {\em longitude}) if it has a slope $1/0$ (resp. $0/1$) and these are denoted by $m$ and $l$ respectively.

\begin{figure}[htbp]
	\begin{center}
	\includegraphics[trim=0mm 0mm 0mm 0mm, width=.6\linewidth]{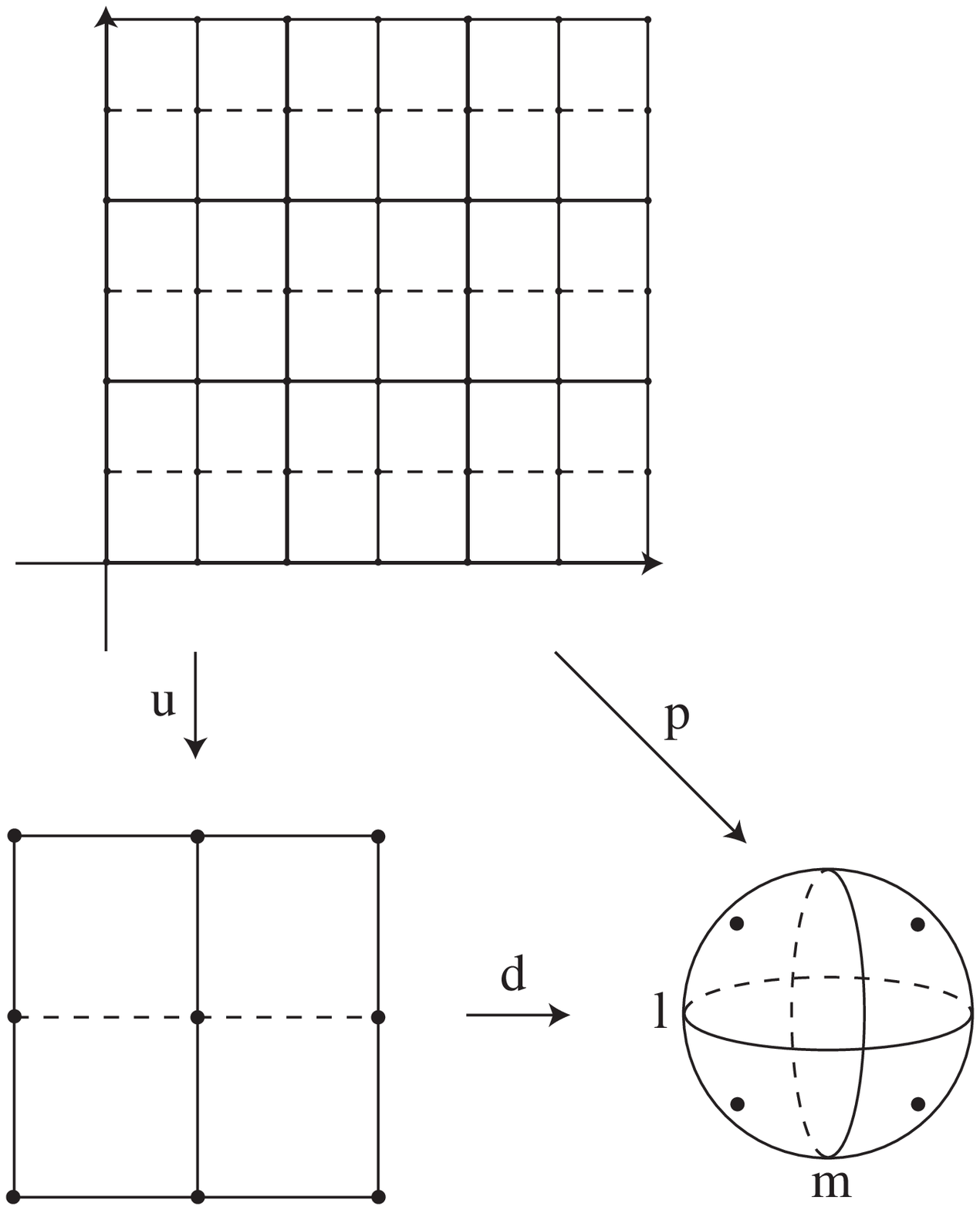}
	\end{center}
	\caption{branched covering map $p:\Bbb{R}^2\to \partial B$}
	\label{covering}
\end{figure}

\begin{lemma}[{\cite[Lemma 3.4]{MO}}]\label{disk}
Let $(B,T)$ be a rational tangle and $F$ an m-essential surface in $B-T$.
Then $F$ is a disk which separates two strings of $T$ in $B$.
\end{lemma}

By Lemma \ref{disk}, a rational tangle $(B,T)$ contains a unique separating disk $D$.
We define the {\em slope} of $(B,T)$ as the slope of $\partial D$.

\subsection{Rotation, reflection, summation and multiplication for tangles}

\begin{figure}[htbp]
	\begin{center}
	\includegraphics[trim=0mm 0mm 0mm 0mm, width=.6\linewidth]{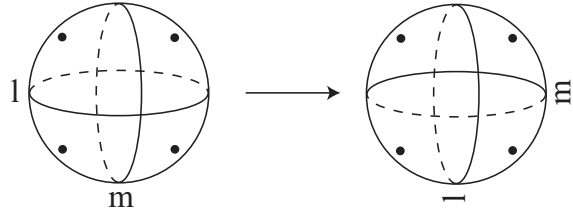}
	\end{center}
	\caption{rotation of a 2-string tangle}
	\label{rotation}
\end{figure}

Let $(B,T)$ be a 2-string tangle fixing four points of $\partial T$.
The {\em rotation} $(B,T^*)$ of $(B,T)$ is obtained by rotating $(B,T)$ counterclockwise by $90^\circ$.
A rotation exchanges the meridian for the longitude.
See Figure \ref{rotation}.

\begin{figure}[htbp]
	\begin{center}
	\includegraphics[trim=0mm 0mm 0mm 0mm, width=.5\linewidth]{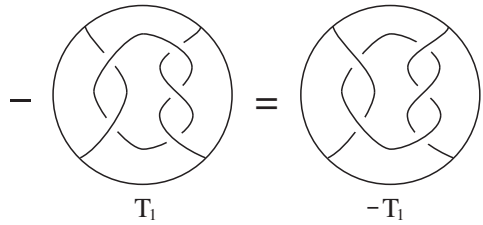}
	\end{center}
	\caption{reflection of a 2-string tangle}
	\label{reflection}
\end{figure}

The {\em reflection} $(B,-T)$ of $(B,T)$ is obtained by reflecting $(B,T)$ by a plane perpendicular to the axis of rotation.
A reflection does not exchange both of the meridian and longitude as sets.
See Figure \ref{reflection}.

\begin{figure}[htbp]
	\begin{center}
	\includegraphics[trim=0mm 0mm 0mm 0mm, width=.9\linewidth]{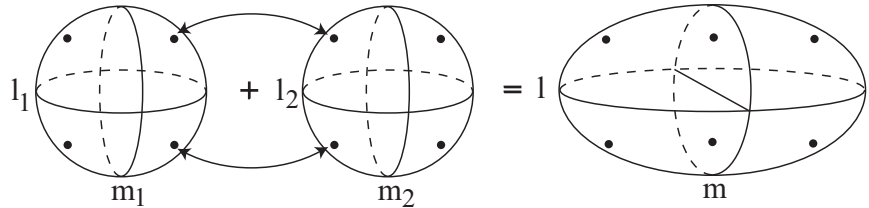}
	\end{center}
	\caption{summation of two 2-string tangles}
	\label{tangle_sum}
\end{figure}

Let $(B_1,T_1)$ and $(B_2,T_2)$ be two tangles with meridians $m_1$, $m_2$ and longitudes $l_1$, $l_2$ respectively.
We construct a new tangle $(B,T)$ from $(B_1,T_1)$ and $(B_2,T_2)$ by gluing the east side disk $D_e$ in $\partial B_1$ and the west side disk $D_w$ in $\partial B_2$ so that $m_1=m_2$, $\partial T_1\cap D_e=\partial T_2\cap D_w$ and $l_1\cap D_e=l_2\cap D_w$.
See Figure \ref{tangle_sum}.
We call this summation a {\em tangle sum} of $(B_1,T_1)$ and $(B_2,T_2)$.
A tangle sum of two tangles $(B_1,T_1)$ and $(B_2,T_2)$ is {\em non-trivial} if neither $(B_i,T_i)$ is a rational tangle of slope $0$ or $\infty$.
An {\em algebraic tangle} is obtained inductively from rational tangles by non-trivial tangle sums and rotations.

Next we define the product of two 2-string tangles $(B_1,T_1)$ and $(B_2,T_2)$.
We regard each crossing of $T_1$ as a rational tangle of slope $1$ or $-1$, and replace each crossing of $T_1$ with $(B_2,T_2)$ or $(B_2,-T_2)$ if the slope of the crossing is $1$ or $-1$ respectively.
Then we obtain the {\em multiplication} $(B,T_1*T_2)$ of two tangles $(B_1,T_1)$ and $(B_2,T_2)$.
We remark that the multiplication depends on the choice of tangle diagrams.
Also in general, $T_1*T_2\ne T_2*T_1$, however two slopes of them coincide (Theorem \ref{homomorphism}).

\begin{figure}[htbp]
	\begin{center}
	\includegraphics[trim=0mm 0mm 0mm 0mm, width=.8\linewidth]{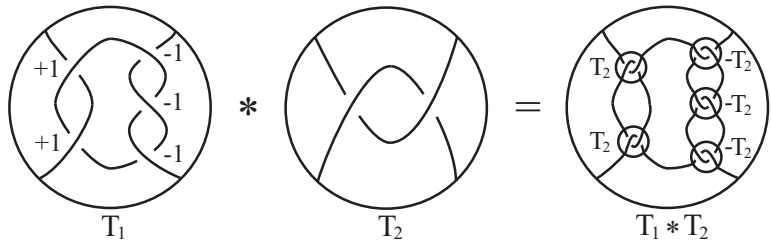}
	\end{center}
	\caption{multiplication of two 2-string tangles}
	\label{multiplication}
\end{figure}

\subsection{The slope of algebraic tangles}

Let $(B,T)$ be a tangle sum of $(B_1,T_1)$ and $(B_2,T_2)$.
Let $F_i$ be a surface in $B_i-T_i$ whose boundary is empty or separates $\partial T_i$ into two pairs of two points.
We say that a surface $F$ in $B-T$ is a {\em non-trivial sum} of $F_1$ and $F_2$ if $F=F_1\cup F_2$, and both of $\partial F_i\cap D$ and $\partial F_i\cap D_i$ are essential arcs or essential loops parallel to $\partial D$ in $D-\partial T_i$ and $D_i-\partial T_i$ respectively, where $D$ is the gluing disk and $D_i=\partial B_i-\rm{int}D$.
The next lemma asserts that a non-trivial sum preserves the m-essentiality of surfaces.
We say that an algebraic tangle is {\em closed} if it is a sum of two tangles both of which have slope $1/0$.

\begin{lemma}\label{sum}
Let $(B,T)$ be a non-trivial sum of two tangles $(B_i,T_i)$.
Then, for any m-essential surface $F$ in $B-T$, there exists an m-essential surface $F_i$ in $B_i-T_i$ such that $F$ is a non-trivial sum of $F_1$ and $F_2$.
Conversely, for any m-essential surface $F_i$ in $B_i-T_i$, a surface obtained from adequate parallelisms of $F_i$ by a non-trivial sum is m-essential in $B-T$.
Moreover, there exists a closed m-essential surface in $B-T$ if and only if there exists a closed algebraic sub-tangle in $(B,T)$.
\end{lemma}

\begin{proof}
Let $F$ be an m-essential surface in $B-T$ and $D$ a disk in $B$ which separates $(B,T)$ into $(B_1,T_1)$ and $(B_2,T_2)$.
Suppose that $F$ intersects $D$ transversely and the intersection is minimal.
Then both of $F_1=F\cap B_1$ and $F_2=F\cap B_2$ are m-essential in $B_1$ and $B_2$ respectively.

Conversely, let $F_i$ be an m-essential surface in $B_i-T_i$ whose boundary runs in $D-T$ as essential arcs or loops $L_i$.
We make a parallelism of $F_i$ so that the number of essential arcs or loops is equal to $l.c.m(|L_1|,|L_2|)$.
Then we obtain an m-essential surface $F$ by a non-trivial sum of a parallelism of $F_1$ and $F_2$.

Moreover, if there exists a closed m-essential surface $F$ in $B-T$, then there exists a minimal algebraic sub-tangle in $(B',T')$ which contains $F$.
We note that $(B',T')$ is not a rational tangle and put $(B',T')=(B_1',T_1')+(B_2',T_2')$.
Then both of $(B_1',T_1')$ and $(B_2',T_2')$ have slope $1/0$, hence $(B',T')$ is closed.
Conversely, if there exists a closed algebraic sub-tangle $(B',T')$ in $(B,T)$, then $B'-T'$ contains a closed m-essential surface $F$ and hence $B-T$ contains $F$.

We omit an elementary cut-and-paste argument required in this lemma.
\end{proof}

By Lemma \ref{sum}, we obtain that the set of all m-essential surfaces in $B-T$ is equal to the set of non-trivial sums of all m-essential surfaces in $B_i-T_i$.
Next theorem extends Lemma \ref{disk}.

\begin{theorem}\label{separating}
Let $(B,T)$ be an algebraic tangle.
Then, any m-essential surface with boundary in $B-T$ separates the components of $T$, and all boundary slopes of m-essential surfaces with boundary in $B-T$ are unique.
\end{theorem}

\begin{proof}
We prove Theorem \ref{separating} by induction on the length of algebraic tangles.
If the length of $(B,T)$ is equal to $1$, then it follows from Lemma \ref{disk}.
Suppose that Theorem \ref{separating} holds for algebraic tangles whose length is less than or equal to $n$, and let $(B,T)$ be an algebraic tangle with length $n+1$.
Then $(B,T)$ is a tangle sum of two algebraic tangles $(B_1,T_1)$ and $(B_2,T_2)$ of slopes $p_1/q_1$ and $p_2/q_2$ whose length is less than or equal to $n$.
By Lemma \ref{sum}, any m-essential surface $F$ in $B-T$ is obtained from parallelisms $m_1F_1$ and $m_2F_2$ of two m-essential surfaces $F_1$ and $F_2$ in $B_1-T_1$ and $B_2-T_2$ by a non-trivial sum.
Here, we note that at least one of $m_1$ and $m_2$ is odd since $F$ is connected.
Suppose without loss of generality that $m_1$ is odd.
Let $X$ and $Y$ be two components that is obtained from $B$ by cutting along $F$.
Then both of $X$ and $Y$ contain some components of $T$ since $F_1$ separates the components of $T_1$ and $m_1$ is odd.
Hence $F$ separates the components of $T$.


Next, suppose that the boundary slopes of m-essential surfaces in $B_i-T_i$ are unique, and let $p_i/q_i$ be the boundary slope of m-essential surfaces in $B_i-T_i$.
If either $q_1$ or $q_2$ is equal to $0$, then by Lemma \ref{sum} any m-essential surface in $B-T$ has the boundary slope $\infty$.
Otherwise, the boundary slope of an m-essential surface $F$ obtained by a non-trivial sum of $F_1$ and $F_2$ is $p_1/q_1+p_2/q_2$.
Hence the boundary slopes of m-essential surfaces in $B-T$ are unique.
\end{proof}

By virtue of Theorem \ref{separating}, we can define the {\em slope} of an algebraic tangle $(B,T)$ as the slope of the boundary of an m-essential surface in $B-T$.

\begin{theorem}\label{homomorphism}
Let $\phi$ be a map from the set of algebraic tangles to the set of rational numbers which maps an algebraic tangle $T$ to the slope of $T$.
Then $\phi$ satisfies the following.
\begin{enumerate}
\item $\phi(T_1+T_2)=\phi(T_1)+\phi(T_2)$ and $\phi(R(0/1))=0$
\item $\phi(T_1*T_2)=\phi(T_1)\phi(T_2)$ and $\phi(R(1/1))=1$
\item $\phi(-T_1)=-\phi(T_1)$
\item $\phi(T_1)\phi(T_1^*)=-1$
\end{enumerate}
where $R(p/q)$ denotes the rational tangle of slope $p/q$.
\end{theorem}

\begin{proof}

(1)  $\phi(T_1+T_2)=\phi(T_1)+\phi(T_2)$ follows the proof of Lemma \ref{sum} and Theorem \ref{separating}.

(3) $\phi(-T_1)=-\phi(T_1)$ holds since a reflection changes the boundary slope  $-1$ times.

(4) $\phi(T_1)\phi(T_1^*)=-1$ follows that a rotation changes the boundary slope reciprocally and times $-1$ by an orthogonal condition.

(2) $\phi(T_1*T_2)=\phi(T_1)\phi(T_2)$ is shown by an induction of the length of $T_1$.
First suppose that $T_1$ is a non-trivial sum of $T_{11}$ and $T_{12}$.
Then,

\begin{align*}
\phi(T_1*T_2) & =\phi((T_{11}+T_{12})*T_2)\\
&=\phi(T_{11}*T_2+T_{12}*T_2)\\
&=\phi(T_{11}*T_2)+\phi(T_{12}*T_2)\\
&=\phi(T_{11})\phi(T_2)+\phi(T_{12})\phi(T_2)\\
&=(\phi(T_{11})+\phi(T_{12}))\phi(T_2)\\
&=\phi(T_{11}+T_{12})\phi(T_2)\\
&=\phi(T_1)\phi(T_2)\\
\end{align*}

Next suppose that $T_1$ is not a non-trivial sum but $T_1^*$ is a non-trivial sum.
Then,

\begin{align*}
\displaystyle \phi(T_1*T_2) &=-\frac{1}{\phi((T_1*T_2)^*)}\\
&=-\frac{1}{\phi(-T_1^**T_2^*)}\\
&=\frac{1}{\phi(T_1^**T_2^*)}\\
&=\frac{1}{\phi(T_1^*)\phi(T_2^*)}\\
&=(-\phi(T_1))(-\phi(T_2))\\
&=\phi(T_1)\phi(T_2)\\
\end{align*}
\end{proof}

\begin{figure}[htbp]
	\begin{center}
	\includegraphics[trim=0mm 0mm 0mm 0mm, width=.6\linewidth]{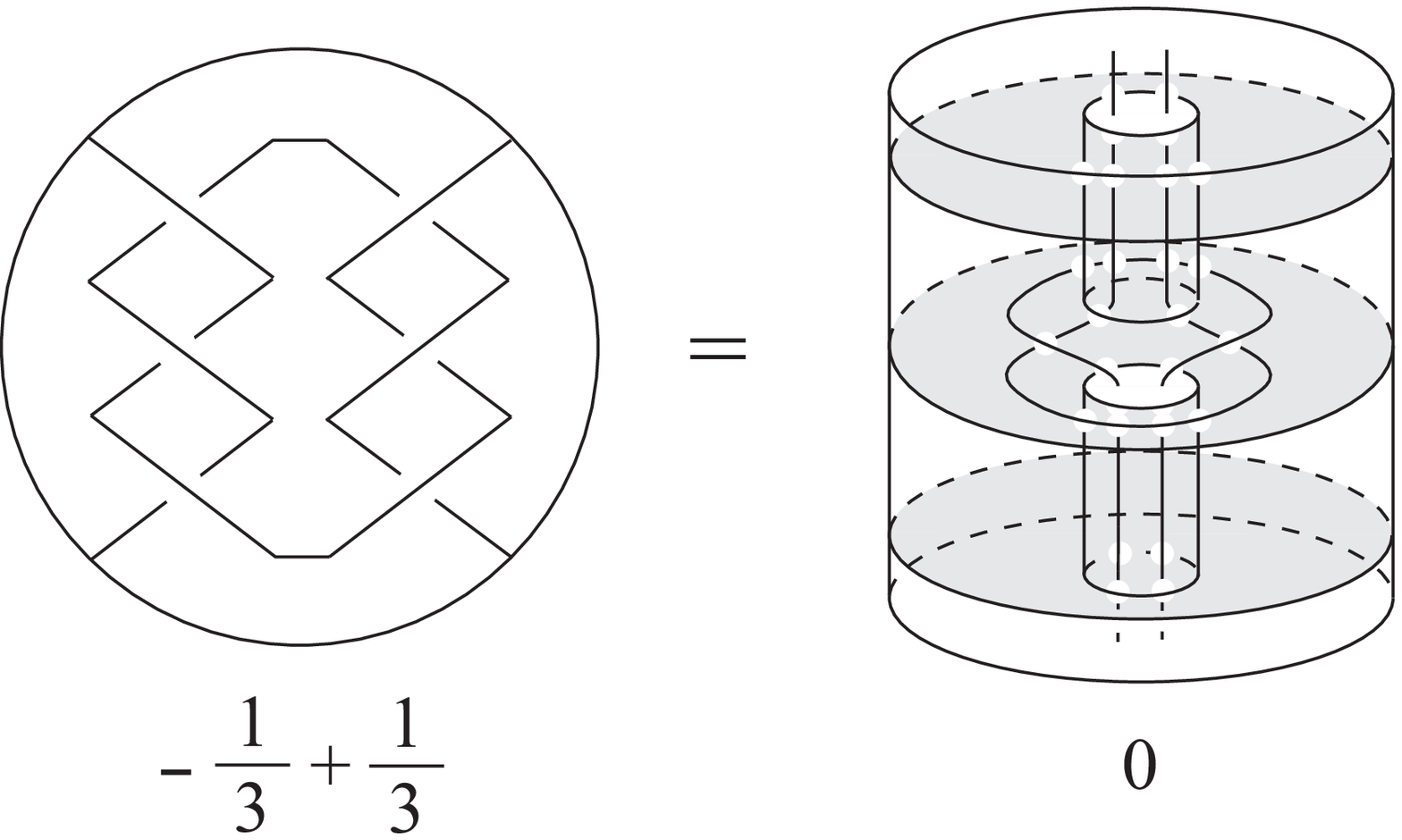}
	\end{center}
	\begin{center}
	\includegraphics[trim=0mm 0mm 0mm 0mm, width=.8\linewidth]{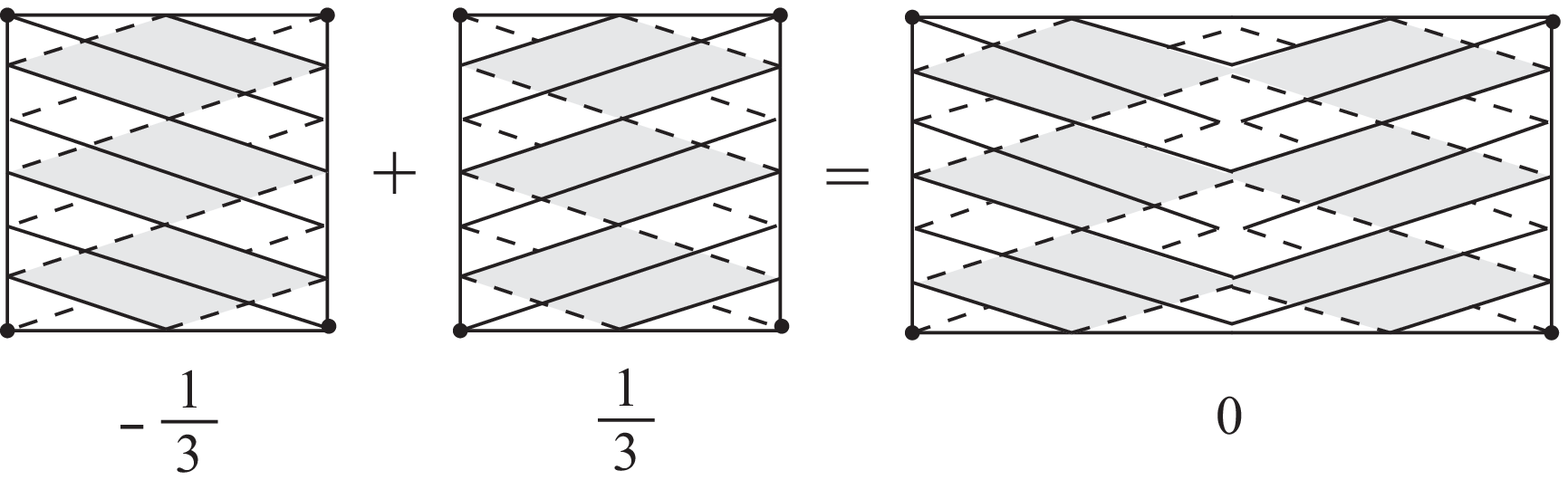}
	\end{center}
	\caption{an algebraic tangle with slope $-1/3+1/3=0$}
	\label{1/3}
\end{figure}

\begin{example}
Figure \ref{1/3} illustrates a non-trivial sum of two rational tangles of slope $-1/3$ and $1/3$.
The resultant algebraic tangle has slope $-1/3+1/3=0$ and contains an m-essential pair of pants which separates two strings.
\end{example}

\begin{figure}[htbp]
	\begin{center}
	\includegraphics[trim=0mm 0mm 0mm 0mm, width=.6\linewidth]{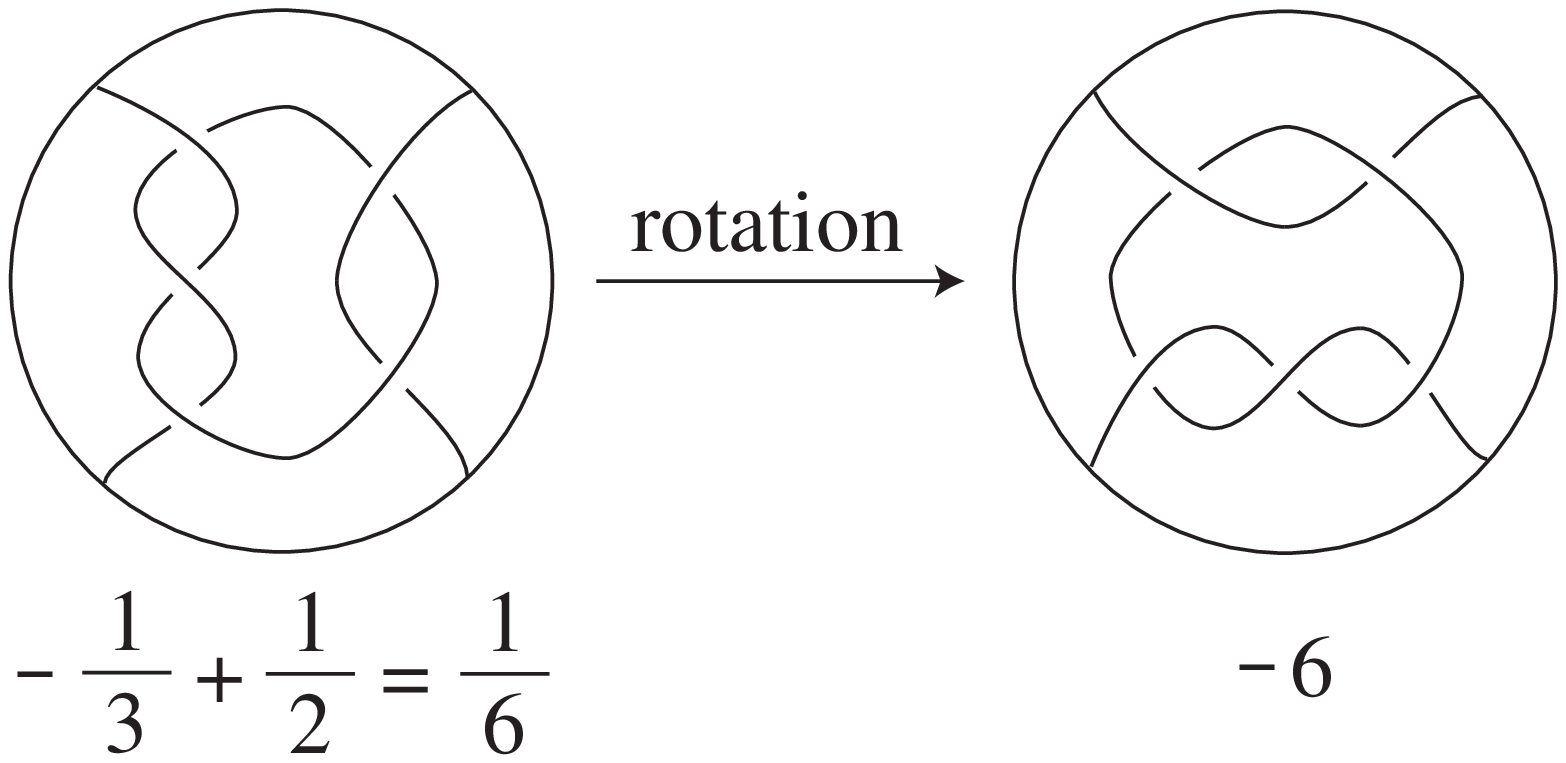}
	\end{center}
	\begin{center}
	\includegraphics[trim=0mm 0mm 0mm 0mm, width=.7\linewidth]{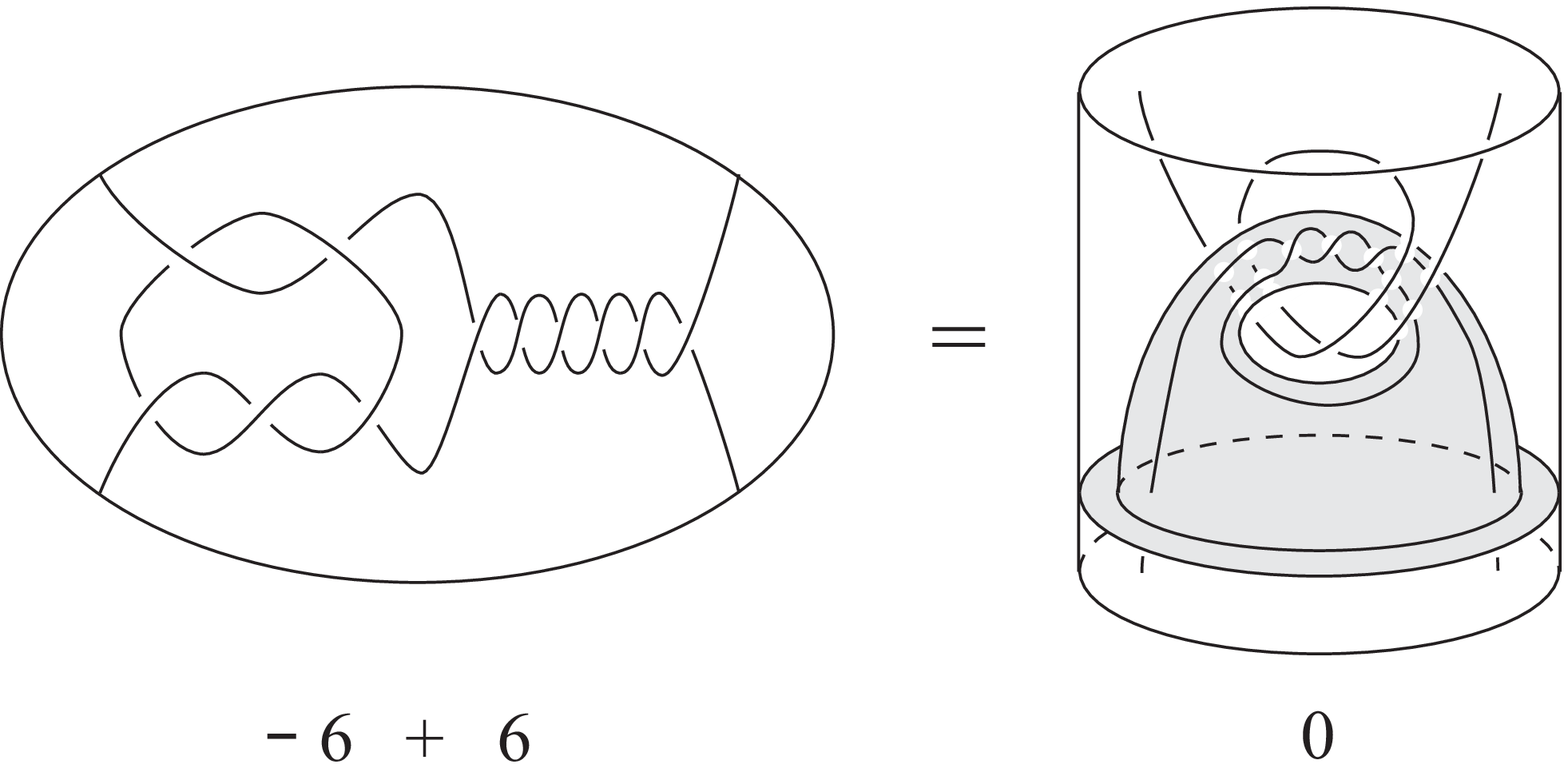}
	\end{center}
	\caption{an algebraic tangle with slope $-1/(1/2+(-1/3))+6=0$}
	\label{1/2-1/3}
\end{figure}

\begin{example}
Figure \ref{1/2-1/3} illustrates a non-trivial sum of three rational tangles of slope $1/2$, $-1/3$ and $6$.
The resultant algebraic tangle has slope $-1/(1/2+(-1/3))+6=0$ and contains an m-essential once punctured torus which separates two strings.
This example is borrowed from $R[1/2,-1/3;-1/6]$ in \cite[Fig. 4.1]{W}.
\end{example}

We say that an algebraic tangle $(B,T)$ of slope $p/q$ is {\em of Type $0/1$} (resp. {\em of Type $1/0$}, {\em of Type $1/1$}) if $p$ is even and $q$ is odd (resp. $p$ is odd and $q$ is even, both of $p$ and $q$ are odd).

\begin{lemma}\label{connection}
An algebraic tangle $(B,T)$ has a loop component if and only if it contains an algebraic sub-tangle which is obtained from two algebraic tangles of Type $1/0$ by a tangle sum.
If an algebraic tangle $(B,T)$ of Type $p/q$ has no loop component, then the connection of two strings $T$ coincides with the rational tangle of slope $p/q$, where $p/q=0/1$, $1/0$ or $1/1$.
\end{lemma}

\begin{proof}
By an elementary calculus on fraction, we have Table \ref{type} on the tangle sums of two algebraic tangles of three Types $0/1$, $1/0$ or $1/1$.
If we do not a tangle sum of two tangles of Type $1/0$, then the connection of two strings $T$ coincides with the rational tangle of slope $p/q$, where $p/q=0/1$, $1/0$ or $1/1$.
Otherwise, the algebraic tangle contains a loop component, and the converse holds.
\end{proof}

\begin{table}[ht]
\noindent\[
\begin{array}{c|c|c|c}
 & $0/1$ & $1/0$ & $1/1$\\
\hline
$0/1$ & $0/1$ & $1/0$ & $1/1$\\
\hline
$1/0$ & $1/0$ & \text{indefinite} & $1/0$\\
\hline
$1/1$ & $1/1$ & $1/0$ & $0/1$
\end{array}
\]
\renewcommand\arraystretch{1.5}
\caption{tangle sums of Types}\label{type}
\end{table}

\subsection{Algebraically alternating knots and links}

We use the next lemma in the proof of Theorem \ref{closed}.

\begin{lemma}[{\cite[Addendum 3.23]{R}}]
Let $(B,T)$ be an algebraic tangle.
\begin{enumerate}
\item $B-T$ contains no m-essential 2-sphere.
\item $B-T$ contains an m-essential disk $D$ if and only if $(B,T)$ is a rational tangle and $D$ is the disk separating the two strings of $T$.
\item $B-T$ contains an m-essential annulus $A$ if and only if $(B,T)=Q_m+(B',T')$ for some $m\ge 1$ and $A$ is a standard annulus in $Q_m$.
\item $B-T$ contains an m-essential torus $F$ if and only if $B-T$ contains a $Q_m$ for some $m\ge 2$ and $F$ is a standard torus in $Q_m$.
\end{enumerate}
\end{lemma}

\begin{proof} (of Theorem \ref{closed})
In this proof, we basically follow Menasco's argument (\cite{M}).
Let $K$ be an algebraically knot or link.
We assume that $K$ is in a position with respect to the 2-sphere $S^2$ as follows.
\begin{enumerate}
\item each algebraic tangle $(B_i,T_i)$ in $(S^3,K)$ intersects $S^2$ in an equatorial disk.
\item the rest of $K$ except all $(B_i,T_i)$ is entirely contained in $S^2$.
\end{enumerate}
We put $S^3=B^3_+\cup_{S^2}B^3_-$, $B_{\pm}=B^3_{\pm}-\rm{int}\bigcup_iB_i$, and $S_{\pm}=\partial B_{\pm}$.

Next, suppose that there exists an m-essential closed surface in $S^3-K$, and let $F$ be an m-essential closed surface in $S^3-K$ such that $\sum_i |F\cap \partial B_i|$ is minimal among all m-essential closed surfaces.
It follows that each component of $F\cap B_i$ is m-essential in $B_i-T_i$ and by Theorem \ref{separating}, it separates $T_i$ in $B_i$ and has the boundary slope $\phi(T_i)$. See Figure \ref{2over3}.
If $F$ is contained in an algebraic tangle $(B_i,T_i)$ in $(S^3,K)$, then by Lemma \ref{sum}, $(B_i,T_i)$ contains a closed algebraic sub-tangle.

\begin{figure}[htbp]
	\begin{center}
	\includegraphics[trim=0mm 0mm 0mm 0mm, width=.3\linewidth]{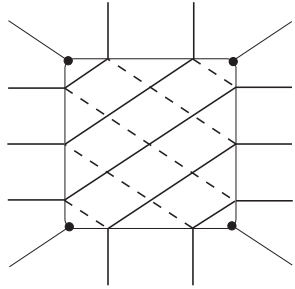}
	\end{center}
	\caption{an intersection of $F$ and the boundary of an algebraic tangle of slope $2/3$}
	\label{2over3}
\end{figure}


By the incompressibility of $F$ in $S^3-K$, we may assume that there is no loop of $F\cap S_{\pm}$ which is entirely contained in $S_+\cap S_-$, and that $F\cap B_{\pm}$ consists of disks.


\begin{claim}\label{once}
There is no loop of $F\cap S_{\pm}$ which runs same algebraic tangle more than once.
\end{claim}

\begin{proof}
Let $\mathcal{L}$ be the set of loops of $F\cap S_+$ which runs same algebraic tangle more than once.
Then, a loop $l$ in $\mathcal{L}$ which is innermost on $S_+$ bounds a disk in $F\cap B_+$.
This shows that $F\cap (S^3-\bigcup_i B_i)$ is boundary compressible in $(S^3-\bigcup_i B_i)-K$, and hence $F$ is compressible or meridionally compressible in $S^3-K$.
\end{proof}

We fix a checkerboard coloring of $(S^2-\bigcup_i B_i)-K$ so that it looks as Figure \ref{coloring} on the basic knot or link $K_0$.
By the alternating property of the basic knot or link $K_0$, a loop of $F\cap S_{\pm}$ satisfies the semi-alternating property;

\begin{enumerate}
\item[(*)] A loop of $F\cap S_+$ goes through over $\partial B_i$ looking $\partial B_i$ in the right side when it runs from a black region to a white region.
\end{enumerate}

\begin{figure}[htbp]
	\begin{center}
	\begin{tabular}{cc}
	\includegraphics[trim=0mm 0mm 0mm 0mm, width=.3\linewidth]{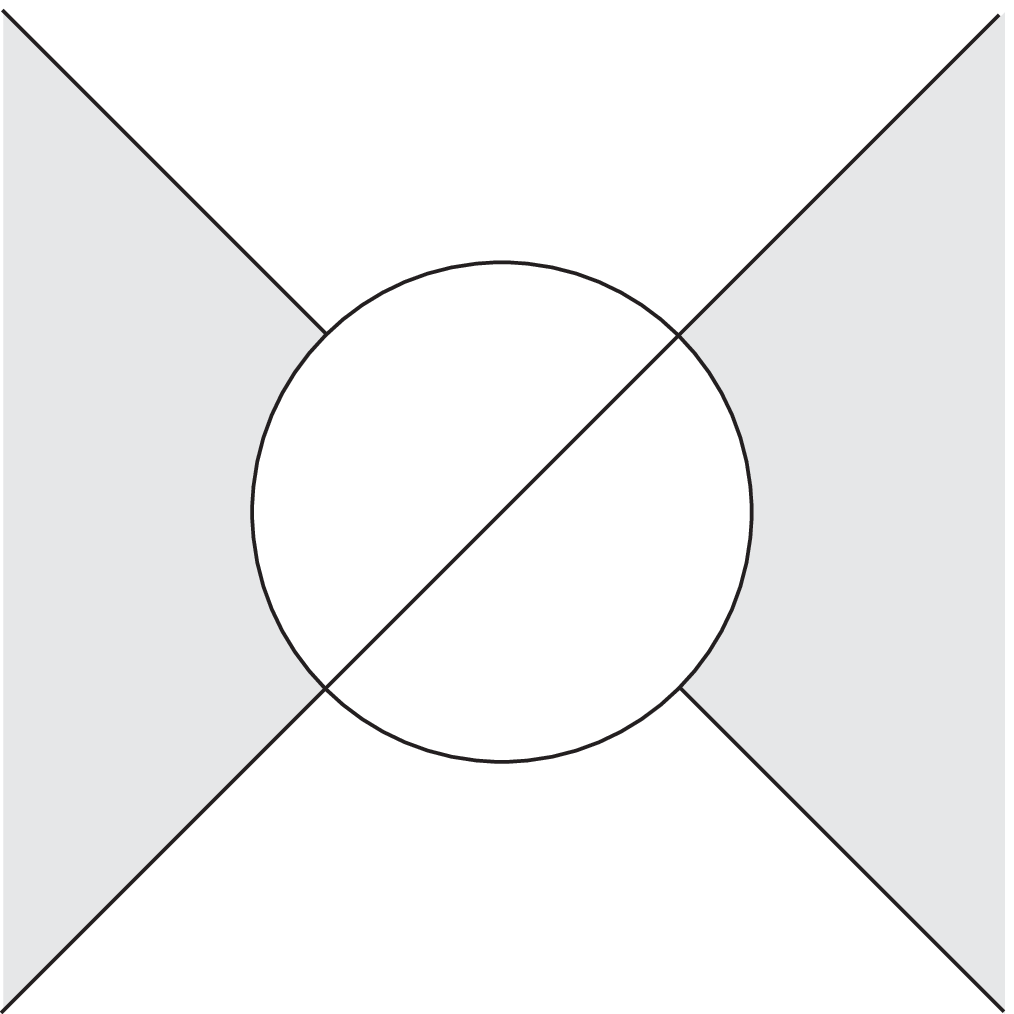}&
	\includegraphics[trim=0mm 0mm 0mm 0mm, width=.3\linewidth]{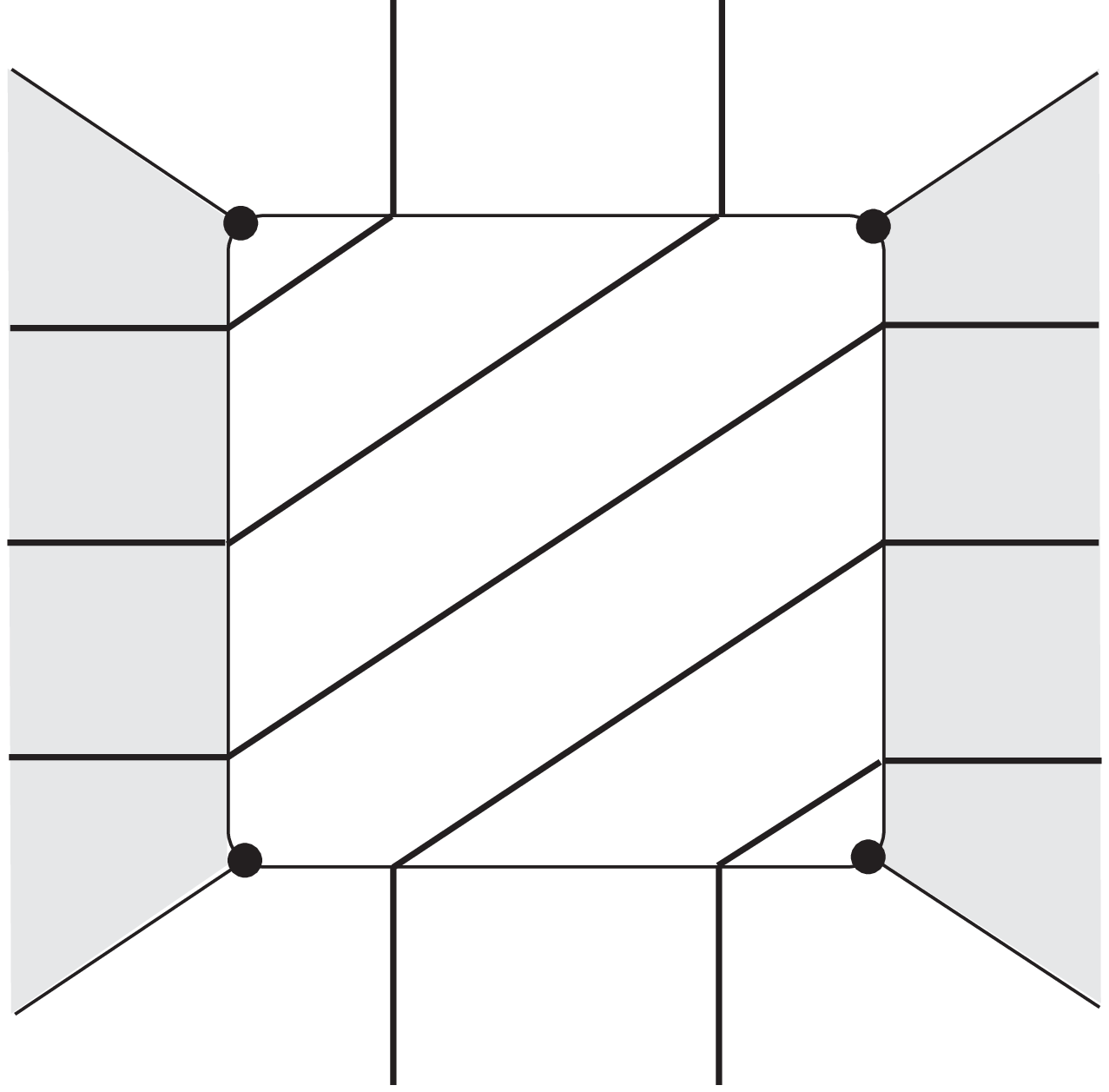}
	\end{tabular}
	\end{center}
	\caption{checkerboard coloring of $K_0$, and loops of $F\cap S_+$ on an algebraic tangle of slope $2/3$}
	\label{coloring}
\end{figure}

\begin{claim}\label{zero}
There is no loop of $F\cap S_+$ which runs through different colors.
\end{claim}

\begin{proof}
Let $\mathcal{L}$ be the set of loops of $F\cap S_+$ which runs through different colors.
Then, a loop in $\mathcal{L}$ which is innermost on $S_+$ never satisfy Claim \ref{once}.
\end{proof}

By Claim \ref{zero}, any loop of $F\cap S_{\pm}$ runs algebraic tangles of slope $0$ or $\infty$.
Hence, $\tilde{K_0}$ is split.
Conversely, if $\tilde{K_0}$ is split, then we can construct an m-essential closed surface in $S^3-K$ which is a union of m-essential surfaces in $\bigcup_i (B_i-T_i)$ and some planar surfaces in the outside of $B_i$.
If there exists an algebraic tangle $(B_i,T_i)$ in $\tilde{K}$ which contains a closed algebraic sub-tangle, then by Lemma \ref{sum}, there exists a closed m-essential surface in $B_i-T_i$ and hence in $S^3-K$.
This proves Theorem \ref{closed} (1).


(2) Suppose that $F$ is an m-essential 2-sphere.
We first note that $F$ is not contained in $\bigcup_i B_i$ by \cite[Addendum 3.23 (1)]{R}.
Then an innermost disk $D$ in $F$ with respect to $F\cap \bigcup_i \partial B_i$ is in the outside of $\bigcup_i B_i$.
By Claim \ref{zero}, $\partial D$ has a slope $0$ or $\infty$ on a $\partial B_i$, and hence the algebraic tangle $(B_i,T_i)$ is a genus 0 cut tangle.

Converserly, if there exists a genus 0 cut tangle $(B_i,T_i)$, then we can construct an m-essential 2-sphere which is a union of an m-essential planar surface in $B_i-T_i$ and some disks in the outside of $B_i$.

(3) Suppose that there exists no genus 0 cut tangle and $F$ is an m-essential torus.
We first note that if $F$ is contained in $\bigcup_i B_i$, then by \cite[Addendum 3.23 (4)]{R}, an algebraic tangle $(B_i,T_i)$ containing $F$ contains $Q_2$.

If there exists a loop of $F\cap \bigcup_i \partial B_i$ which is inessential in $F$, then by \cite[Addendum 3.23 (2)]{R}, an innermost disk $D$ on $F$ is in the outside of $\bigcup_i B_i$.
By Claim \ref{zero}, $\partial D$ has a slope $0$ or $\infty$ on a $\partial B_i$, and hence the algebraic tangle $(B_i,T_i)$ is a genus 1 cut tangle.
Otherwise, since the basic diagram $\tilde{K_0}$ is connected, each component of $F\cap \bigcup_i \partial B_i$ is an essential loop in $F$, and hence every loops are parallel in $F$.
Let $A_1$ be an annulus of $F\cap \bigcup_i B_i$.
By \cite[Addendum 3.23 (3)]{R}, the algebraic tangle $B_i$ containing $A_1$ is a tangle sum of $Q_m$ $(m\ge 1)$ and an algebraic tangle, and $A_1$ is a standard annulus in $Q_m$.
Let $A_2$ be the next annulus of $A_1$ in $F$.
By Claim \ref{once}, $A_2$ connects $B_i$ and other algebraic tangle $B_j$.
Then the next annulus $A_3$ of $A_2$ in $F$ is an essential annulus in $(B_j,T_j)$, and by \cite[Addendum 3.23 (3)]{R}, the algebraic tangle $(B_j,T_j)$ containing $A_3$ is a tangle sum of $Q_{m'}$ $(m'\ge 1)$ and an algebraic tangle, and $A_3$ is a standard annulus in $Q_{m'}$.
A loop of $Q_{m'}$ is parallel to a loop of $Q_m$ along $A_2$ and hence $(S^3,K)$ contains $Q_2$.

Conversely, if there exists a genus 1 cut tangle $(B_i,T_i)$, then we can construct an m-essential torus which is a union of an m-essential genus 1 surface in $B_i-T_i$ and some disks in the outside of $B_i$.
\end{proof}

\bigskip


\bibliographystyle{amsplain}

\end{document}